\documentclass[leqno,11pt,a4paper]{amsart}
\usepackage{graphicx}
\usepackage[usenames,dvipsnames]{color}
\usepackage{hyperref}
\usepackage{booktabs}
\usepackage{array,longtable}
\newtheorem{thm}{Theorem}[section]
\newtheorem{eg}[thm]{Example}
\newtheorem{remark}[thm]{Remark}
\numberwithin{equation}{section}
\newenvironment{acknowledge}{\bigskip\noindent\textbf{Acknowledgments.}}{}
\long\def\blankfootnotetext#1{\begingroup\def\thefootnote{\fnsymbol{footnote}}\footnotetext{#1}\endgroup}
\newcommand{\Z}{\mathbb{Z}}
\newcommand{\Q}{\mathbb{Q}}
\newcommand{\F}{\mathbb{F}}
\newcommand{\cP}{\mathcal{P}}
\newcommand{\cS}{\mathcal{S}}
\newcommand{\GL}{\mathrm{GL}}
\newcommand{\abs}[1]{\left\vert{#1}\right\vert}
\newcommand{\V}[1]{\mathrm{vert}\!\left({#1}\right)}
\newcommand{\conv}[1]{\mathrm{conv}\!\left({#1}\right)}
\newcommand{\Aut}[1]{\mathrm{Aut}\!\left({#1}\right)}
\newcommand{\AffAut}[1]{\mathrm{AffAut}\!\left({#1}\right)}
\newcommand{\magma}{{\sc Magma}}
\graphicspath{{images/}}
\begin{document}
\author[G.~Brown]{Gavin Brown}
\address{Department of Mathematical Sciences\\Loughborough University\\Loughborough\\LE11 3TU\\United Kingdom}
\email{G.D.Brown@lboro.ac.uk}
\author[A.~M.~Kasprzyk]{Alexander M.~Kasprzyk}
\address{Department of Mathematics\\Imperial College London\\London\\SW7 2AZ\\United Kingdom}
\email{a.m.kasprzyk@imperial.ac.uk}
\blankfootnotetext{2010 \emph{Mathematics Subject Classification}: 14G50 (Primary); 52B20, 14M25 (Secondary).}
\title{Seven new champion linear codes}
\begin{abstract}
We exhibit seven linear codes exceeding the current best known minimum distance $d$ for their dimension $k$ and block length $n$. Each code is defined over $\F_8$, and their invariants $[n,k,d]$ are given by $[49,13,27]$, $[49,14,26]$, $[49,16,24]$, $[49,17,23]$, $[49,19,21]$, $[49,25,16]$ and $[49,26,15]$. Our method includes an exhaustive search of all monomial evaluation codes generated by points in the $[0,5]\times [0,5]$ lattice square.
\end{abstract}
\maketitle
\section{Introduction}\label{sec:introduction}
The basic invariants of a linear code $C$ over the finite field $\F_q$ are its \emph{dimension} $k$ and its \emph{block length} $n$. The code is the image in $\F_q^n$ of a $k\times n$ matrix $M$ over $\F_q$.
The \emph{minimum distance} $d$ of $C$ is the smallest Hamming weight of any nonzero linear combination of the rows of $M$. These invariants are conventionally recorded as a triple $[n,k,d]$. One usually seeks codes with minimum distance $d$ as large as possible for a given block length $n$ and dimension $k$. There are theoretical upper bounds for the minimum distance, although in many cases the largest known example falls short of these bounds. Grassl~\cite{grassl} catalogues the best known linear codes in this sense, and we refer to any code with larger minimum distance as a \emph{champion code}\footnote{The profile $[36,19,12]$ over $\F_7$ discovered in~\cite{BK12} is not currently recorded in~\cite{grassl}, but we will~\emph{not} regard any code matching these invariants as a champion.}.

In this paper we find seven new champion codes defined over $\F_8$ by considering the class of \emph{generalised toric codes} introduced by Little~\cite{Lit11} (see~\S\ref{sec:generalised}).

\begin{thm}\label{thm:exhaustive_search_small_squares}
There are precisely five $[n,k,d]$ profiles of champion generalised toric codes over $\F_8$ generated by collections of points in a $[0,5]\times [0,5]$ square. There are (at least) another two  champion codes arising from a $[0,6]\times [0,6]$ square. These profiles are listed in Table~\ref{tab:champs}.
\end{thm}

The codes we consider correspond to subsets of lattice points in the $[0,m]\times [0,m]$ square for small values of $m$. In principle, our approach is to enumerate all such sets of points, construct the corresponding code over various small fields and then compute its invariants. In~\S\ref{sec:algorithm} we describe an  algorithm to enumerate all such sets of points up to affine equivalence. To prove the theorem, we implement and run this algorithm to completion for $m=5$, and a partial search when $m=6$. Since our algorithm is exhaustive, we can say more:

\begin{thm}\label{thm:exhaustive_search_small_fields}
There are no champion generalised toric codes over $\F_q$ when  $q\le 7$.
\end{thm}

It is worth noting that the new Bound~A of~\cite{GBS} does not improve the existing theoretical upper bounds for minimum distance in the range of the champion codes we discover, so we cannot say that the codes of Theorem~\ref{thm:exhaustive_search_small_squares} are the best possible amongst all linear codes.

\begin{table}[htbp]
\centering
\begin{tabular}{ccccp{4.4in}} \toprule
$n$&$k$&$d$&Min\,$m$&\multicolumn{1}{c}{Example points}\\ 
\cmidrule(lr){1-3} \cmidrule(lr){4-4} \cmidrule(lr){5-5}
$49$&$13$&$27$&$4$&
\small$(0,2),$ $(0,3),$ $(0,4),$ $(1,0),$ $(1,2),$ $(2,1),$ $(2,2),$ $(2,4),$ $(3,0),$ $(3,2),$ $(4,1),$ $(4,3),$ $(4,4)$\\
$49$&$14$&$26$&$5$&
\small$(0,1)$ $(0,4),$ $(1,3),$ $(2,3),$ $(2,4),$ $(3,1),$ $(3,3),$ $(3,5),$ $(4,0),$ $(4,2),$ $(4,4),$ $(4,5),$ $(5,2),$ $(5,5)$\\
$49$&$16$&$24$&$5$&
\small$(0,1),$ $(0,2),$ $(0,4),$ $(0,5),$ $(1,1),$ $(1,4),$ $(1,5),$ $(2,0),$ $(2,1),$ $(2,4),$ $(3,2),$ $(3,3),$ $(3,5),$ $(4,2),$ $(5,0),$ $(5,5)$\\
$49$&$17$&$23$&$5$&
\small$(0,1),$ $(0,3),$ $(0,5),$ $(1,1),$ $(1,2),$ $(1,4),$ $(2,2),$ $(2,3),$ $(2,4),$ $(3,0),$ $(3,1),$ $(3,4),$ $(3,5),$ $(4,0),$ $(4,2),$ $(5,1),$ $(5,5)$\\
$49$&$19$&$21$&$5$&
\small$(0,0),$ $(0,2),$ $(0,4),$ $(0,5),$ $(1,0),$ $(1,1),$ $(1,5),$ $(2,1),$ $(2,2),$ $(2,5),$ $(3,1),$ $(3,3),$ $(3,4),$ $(4,0),$ $(4,5),$ $(5,0),$ $(5,2),$ $(5,4),$ $(5,5)$\\
$49$&$25$&$16$&$6$&
\small$(0,4),$ $(0,5),$ $(0,6),$ $(1,2),$ $(1,3),$ $(1,4),$ $(1,5),$ $(1,6),$ $(2,1),$ $(2,2),$ $(2,3),$ $(2,5),$ $(2,6),$ $(3,0),$ $(3,2),$ $(3,4),$ $(3,6),$ $(4,2),$ $(4,3),$ $(4,4),$ $(4,5),$ $(5,1),$ $(5,5),$ $(6,2),$ $(6,3)$\\
$49$&$26$&$15$&$6$&
\small$(0,2),$ $(0,4),$ $(0,5),$ $(1,1),$ $(1,2),$ $(1,3),$ $(1,4),$ $(1,6),$ $(2,2),$ $(2,3),$ $(2,5),$ $(3,0),$ $(3,1),$ $(3,3),$ $(3,5),$ $(3,6),$ $(4,3),$ $(4,4),$ $(5,0),$ $(5,2),$ $(5,3),$ $(5,4),$ $(5,5),$ $(6,1),$ $(6,2),$ $(6,4)$\\
\bottomrule \\
\end{tabular}
\caption{The $[n,k,d]$ invariants of new champion generalised toric codes over $\F_8$. In each case a single example is given, contained in the smallest possible $[0,m]\times[0,m]$ square, and illustrated in Figure~\ref{fig:champs}.}\label{tab:champs}
\end{table}

\section{Generalised toric codes}\label{sec:generalised}
Recall from~\cite{Han00} that a \emph{toric code} $C$, over a sufficiently large field $\F_q$, is determined by a convex lattice polygon $P\subset\Z^2\otimes\Q$ as follows. Suppose that $P$ lies in a $[0,m]\times [0,m]$ square, where $q\ge m+2$. Then $C$ is given by the image of a $k\times n$ matrix $M$ whose rows are generated by evaluating each lattice point $(a,b)\in P\cap\Z^2$ (regarded as a monomial $x^ay^b$) at each vector of the torus $(\F^*_q)^2$. The dimension of $C$ is the number of rows $k=\abs{P\cap\Z^2}$ of $M$ and the block length is the number of columns $n=(q-1)^2$.

\begin{figure}[htbp]
\centering
\begin{tabular}{cc}
\includegraphics[scale=0.9]{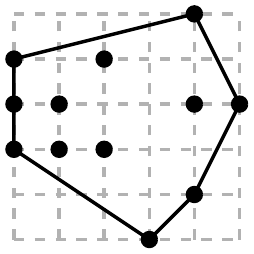}&
\includegraphics[scale=0.9]{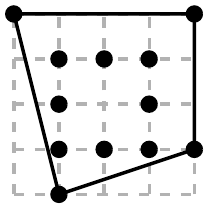}\\
\small(a)&\small(b)
\end{tabular}
\caption{(a) The $[49,12,28]$ code over $\F_8$ described by Little~\cite{Lit11}; (b) One of $448$ non-equivalent ways of generating a $[49,12,28]$ code over $\F_8$ from points in a $[0,4]\times [0,4]$ square.}\label{fig:little}
\end{figure}

A \emph{generalised toric code} is constructed in the same way as a toric code but with the possibility of omitting one or more of the lattice points of $P$; equivalently, one may remove rows from the generating matrix $M$. Conceptually this allows one to delete any particularly short vectors in $C$ that arise as rows of $M$. Although of course this does not necessarily increase the minimum distance, they have recently also been proving fruitful in the search for champions~\cite{Rua09,AHV09,CMO09,Lit11}. Little introduced the study of generalised toric codes, and found a champion $[49,12,28]$ code over $\F_8$ coming from a particular subset of points of a polygon in a $[0,5]\times [0,5]$ square~\cite{Lit11}; this is illustrated in Figure~\ref{fig:little}(a). Systematic attempts to produce examples of champion generalised toric codes have been performed over $\F_4$ and $\F_5$~\cite{AHV09}, and over $\F_7$, $\F_8$ and $\F_9$~\cite{CMO09}.

\begin{figure}[htbp]
\centering
\includegraphics[scale=0.9]{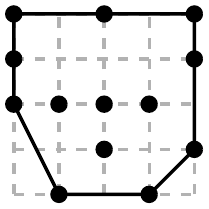}
\\\vspace{7px}
\includegraphics[scale=0.9]{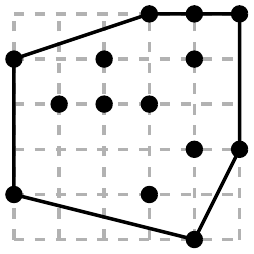}\hspace{7px}
\includegraphics[scale=0.9]{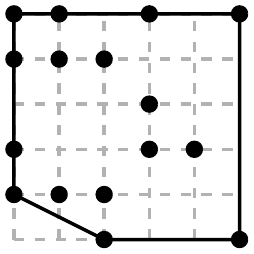}\hspace{7px}
\includegraphics[scale=0.9]{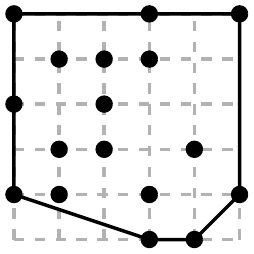}\hspace{7px}
\includegraphics[scale=0.9]{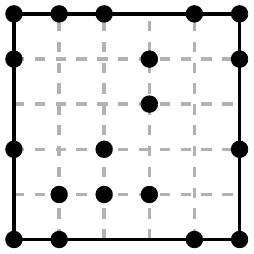}
\\\vspace{7px}
\includegraphics[scale=0.9]{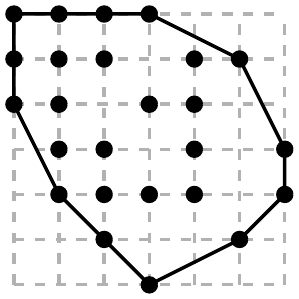}\hspace{7px}
\includegraphics[scale=0.9]{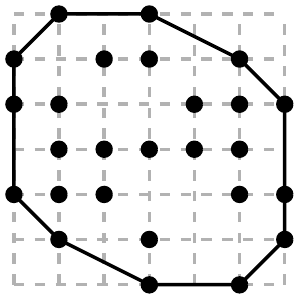}
\caption{Example point configurations for the new champion generalised codes over $\F_8$ listed in Table~\ref{tab:champs}.}\label{fig:champs}
\end{figure}

\subsection*{Small polygons}
In~\cite{BK12} we assembled a comprehensive database of lattice polygons that are contained, up to lattice automorphism and translation, in a $[0,m]\times [0,m]$ square, for $m\le 7$. This database can be interrogated online via the Graded Ring Database~\cite{GRDb}, or from within the computational algebra software {\magma}~\cite{Magma}. In~\cite{BK12} we checked the toric code corresponding to each such polygon, over all prime-powered fields $\F_q$ for $m+2\le q\le 9$, and found a single new champion code, defined over $\F_7$ with invariants $[36,19,12]$. In this paper we consider all \emph{generalised} toric codes associated with polygons that lie in a $[0,m]\times [0,m]$ square, for $m\le 5$ and $m+2\le q\le 9$. We also present partial results when $m=6$; in this case the number of possible generalised toric codes is too large to be searched systematically using current techniques. The invariants listed in Table~\ref{tab:champs} are those of the champion codes we found with this search.

\subsection*{Multiplicity of champions}
\begin{figure}[htbp]
\centering
\includegraphics[scale=0.9]{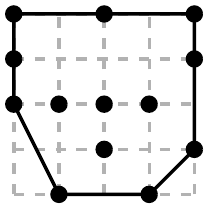}\hspace{7px}
\includegraphics[scale=0.9]{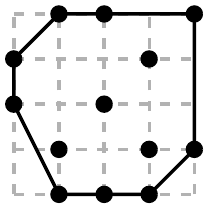}\hspace{7px}
\includegraphics[scale=0.9]{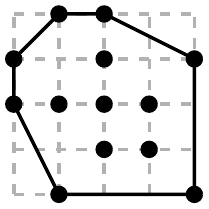}\hspace{7px}
\includegraphics[scale=0.9]{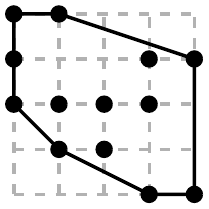}
\caption{The possible choices of points in a $[0,4]\times [0,4]$ square giving a $[49,13,27]$ code over $\F_8$.}\label{fig:nonequivalent_points}
\end{figure}

Champion profiles such as $[49,13,27]$ over $\F_8$ are often achieved in many non-equivalent ways: this case, for example, is realised by four different sets of points in a $[0,4]\times [0,4]$ square, distinct up to lattice automorphism and translation; see Figure~\ref{fig:nonequivalent_points}. Little's champion $[49,12,28]$ code is more striking still: we have $448$ non-equivalent generalised codes in a $[0,4]\times [0,4]$ square that yield these invariants. An example is illustrated in Figure~\ref{fig:little}(b).

\subsection*{Toric codes often achieve best known minimum distance}
We extracted the $[n,k,d]$ profiles of all generalised toric codes that matched, or exceeded, the largest minimum distances available in Grassl's catalogue over $\F_q$ for $q\leq 7$. All such codes are contained in a $[0,5]\times[0,5]$ square, and so our results are complete. The bounds attained are recorded in Appendix~\ref{app:toric_table}; generalised toric codes achieve (or exceed) the current best known minimum distance in $28$ of the $57$ cases.

\section{The algorithm}\label{sec:algorithm}

The proof of Theorems~\ref{thm:exhaustive_search_small_squares} and~\ref{thm:exhaustive_search_small_fields} is by a systematic computer search. The hurdle to overcome is the sheer number of codes and the complexity of computing their minimum distances (often these calculations would take millions of years if the code achieved its apparent minimum distance). We describe an algorithm that carries out the enumeration of subsets of lattice points of a lattice square up to affine equivalence. We have implemented this in {\magma} exactly as described here, making use of the convex polytopes package~\cite{ConvChap}. Our code is available to download from~\cite{Magmacode}.

Let $S\subset\Z^2$ be a collection of points contained in a $[0,m]\times[0,m]$ square generating the generalised toric code $C_S$ over $\F_q$. For any lattice translation $S'=S-u$ of the points $S$, or more generally for any affine linear transformation $s\mapsto (s-u)M$ of $S$, where $u\in\Z^2$ and $M\in\GL_2(\Z)$, the toric codes $C_S$ and $C_{S'}$ are monomially equivalent~\cite[Theorem~4]{LS07}. Therefore it is enough to consider points $S$ up to lattice translation and change of basis. Consider now the lattice polygon $P:=\conv{S}\subset\Z^2\otimes\Q$. Up to equivalence, $P$ can be assumed to be one of the polygons constructed in~\cite{BK12}. This motivates our approach: For each polygon $P$ contained in a $[0,m]\times[0,m]$ square we will generate, up to equivalence, all possible subsets $S$ of points $P\cap\Z^2$ such that $S$ contains the vertices $\V{P}$ of $P$. Insisting that $\V{P}\subseteq S$ is a natural restriction; if this were not so, $Q:=\conv{S}$ is a lattice polygon distinct from $P$ in the same $[0,m]\times[0,m]$ square, and, up to equivalence, that polygon will be considered separately.

\subsection*{Step~1: Compute the affine automorphism group of $P$}
The first step is to compute the affine automorphism group $G:=\AffAut{P}$ of the polygon $P$, where any element $g\in G$ can be written as a combination of elements of $\GL_2(\Z)$ and translations. Embed $P$ at height~$1$ in the lattice $\Z^3 = \Z^2\times\Z$; for example, simply append a coordinate~$1$ to each vertex of $P$. We refer to this embedded image of $P$ also as $P$. Form the cone $\sigma$ with vertex the origin generated by the points of the embedded image of $P$. Let $\Aut{\sigma}$ be the linear automorphism group of $\sigma$. The action of this group on $\sigma$ restricts to a faithful action on $P$, realising the full group of affine lattice automorphisms of~$P$.

\subsection*{Step~2: Extend the action to subsets of points of $P$}
Clearly $v\cdot G\subseteq\V{P}$ for any $v\in\V{P}$. Let $\mathcal{P}:=P\cap\Z^2\setminus\V{P}$ be the set of non-vertex lattice points in~$P$. We fix an order of the points of $\mathcal{P} = \{v_1,\dots,v_k\}$, where $k:=\abs{\mathcal{P}}$; $G$ acts on $\mathcal{P}$ via permutation. Choose a largest orbit $O_1$ of this action on~$\mathcal{P}$. Without loss of generality we assume that $O_1$ permutes the first $k_1\leq k$ elements.

\subsection*{Step~3: Enumerate subsets of a largest orbit up to the action}
From now onwards we regard $G$ as acting on $\{0,1\}^k$, where an element $s\in\{0,1\}^k$ corresponds to a choice of points $S\subseteq P\cap\Z^2$ via the obvious map
$$s=(b_1,\dots,b_k)\mapsto S=\V{P}\cup\{v_i\in\mathcal{P}\mid 1\leq i\leq k, b_i=1\}.$$
Let $\cS_1$ be the set of all sequences in $\{0,1\}^{k_1}$, up to this action. Enumerating $\cS_1$ for small polygons is straightforward: in particular, $\cS_1$ contains a unique sequence of zeros and also a unique sequence containing a single~$1$. Larger numbers $h\ge 2$ of nonzero coefficients depend on the $h$-transitivity of the action: the number is determined by the orbit--stabiliser theorem, but even without applying that the numbers are small enough simply to run through all possibilities rejecting those already seen (up to the action). We only compute the results for at most $\lfloor k_1/2 \rfloor$ nonzero coefficients, since the remaining possibilities are obtained via symmetry by exchanging $0$ and $1$.

\subsection*{Step~4: Extend subsets to all of $\cP$}
For each $s\in\mathcal{S}_1$ we extend $s$ to an element of $\{0,1\}^k$ as follows. Let $G_s:=\{g\in G\mid s\cdot g=s\}\le G$, and let $\mathcal{P}_s:=\{v_i\in\mathcal{P}\mid k_1<i\leq k\}$. By induction we can construct the set of possible $\{0,1\}^{k-k_1}$-sequences $\mathcal{S}_s$, up to the action of $G_s$. Let $s\sqcup s'$ denote the concatenation of $s\in\{0,1\}^{k_1}$ with an element $s'\in\{0,1\}^{k-k_1}$, so that $s\sqcup s'\in\{0,1\}^k$. Then
$$\mathcal{S}=\{s\sqcup s'\mid s\in\mathcal{S}_1, s'\in\mathcal{S}_s\}$$
corresponds to the set of all points in $P\cap\Z^2\setminus\V{P}$, up to the action of $G$.

\begin{figure}[htbp]
\centering
\includegraphics[scale=0.9]{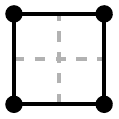}\hspace{7px}
\includegraphics[scale=0.9]{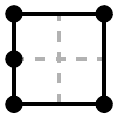}\hspace{7px}
\includegraphics[scale=0.9]{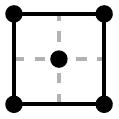}\hspace{7px}
\includegraphics[scale=0.9]{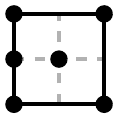}\hspace{7px}
\includegraphics[scale=0.9]{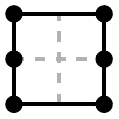}\hspace{7px}
\includegraphics[scale=0.9]{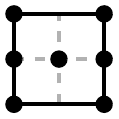}
\\\vspace{7px}
\includegraphics[scale=0.9]{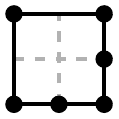}\hspace{7px}
\includegraphics[scale=0.9]{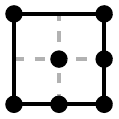}\hspace{7px}
\includegraphics[scale=0.9]{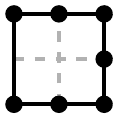}\hspace{7px}
\includegraphics[scale=0.9]{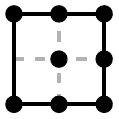}\hspace{7px}
\includegraphics[scale=0.9]{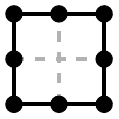}\hspace{7px}
\includegraphics[scale=0.9]{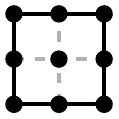}
\caption{The twelve equivalence classes of points $S$ in the $[0,2]\times[0,2]$ square $B$ such that $\V{B}\subseteq S$.}\label{fig:square_points}
\end{figure}

\begin{eg}
\rm
Consider the polygon $B$ equal to the $[0,2]\times[0,2]$ square. This contains $3^2$ points, and hence there exist $2^{3^2}$ subsets of points in $B$. If one insists that $\V{B}\subseteq S$ for any subset $S$ of points, then this is reduced to $2^{3^2-4}=32$ possibilities. Considering subsets $S$ only up to equivalence, we obtain the twelve equivalence classes illustrated in Figure~\ref{fig:square_points}.
\end{eg}

At this point, we have a list of all possible subsets of points of $P$ up to affine automorphisms (and including its vertices). It remains to identify any champion generalised toric codes.

\subsection*{Step~5: Run a trial minimum distance algorithm to exclude most cases}
Fix a prime power $q$ and an integer $0<m<q-1$. Given a set of points $S\subset [0,m]\times [0,m]$, compute the generalised toric code $C = C_S(\F_q)$. (In our case, $S$ is one of the subsets of points of a polygon enumerated in Step~$4$.) Let $M = M_C$ be the generator matrix of $C$. Compute the Hamming length of the shortest nonzero vector that is a linear combination over $\F_q$ of up to four rows of $M$ -- this is $q^4-1$ vectors, which is manageable for the small values of $q$ we consider. If this length is strictly greater than the best known minimum distance for codes of dimension $\abs{S}$ and block length $q-1$, then keep the pair $S,q$ for consideration at the next step; otherwise discard the pair $S,q$, since the corresponding code cannot be a new champion.

\subsection*{Step~6: Compute the minimum distance of successful trial candidates}
In principle this is now the hard part. We employ {\magma}'s {\tt MinimumDistance} function; at heart, it must simply check the Hamming weight of all words in the code, and so can easily be expected to take too long for our purposes. In practice when $m\le 5$, for the hundreds of thousands of codes that pass our four-line trial in Step~$5$, we can always calculate the minimum distance within a few hours, and there are several hundred cases realising champion invariants. So for $m\le 5$ there are no numerical profiles $[n,d,k]$ other than the five listed in Table~\ref{tab:champs} that admit champion generalised codes.

\begin{remark}
\rm
To get an idea of the effect of symmetries, Table~\ref{tab:squares} lists the number of subsets of points $2^{(m+1)^2-4}$ of the $[0,m]\times [0,m]$ square (excluding the vertices, since we handled vertices of polygons separately) and the number of such subsets up to equivalence. Of course there are eight symmetries of the square, and $8\times 265488 = 2123904$; the full group of symmetries of $P$ has been exploited.

\begin{table}[htbp]
\centering
$\begin{array}{llcccc} \toprule
\multicolumn{1}{c}{m}&&1&2&3&4\\ \midrule
\text{\# subsets of points}& \hspace{1ex}& 1 & 32 & 4096 & 2097152\\
\text{\# up to symmetries} && 1 & 12 & 570 &265488\\
\bottomrule \\
\end{array}$
\caption{A comparison of the number of subsets of points in a $[0,m]\times[0,m]$ square, and the number of subsets up to equivalence.}\label{tab:squares}
\end{table}

It is conceivable that there are other `hidden' symmetries that preserve the invariants of the toric code, even if they do not preserve the lattice subset, and the multiplicity of champions may be hinting at this.
\end{remark}

\begin{acknowledge}
Our thanks to Tom Coates for several helpful ideas, to John Cannon for providing {\magma} for use on the Imperial College mathematics cluster, and to Andy Thomas for technical assistance. The second author is supported by EPSRC grant EP/I008128/1.
\end{acknowledge}

\bibliographystyle{amsalpha}
\providecommand{\bysame}{\leavevmode\hbox to3em{\hrulefill}\thinspace}
\providecommand{\MR}{\relax\ifhmode\unskip\space\fi MR }
\providecommand{\MRhref}[2]{%
  \href{http://www.ams.org/mathscinet-getitem?mr=#1}{#2}
}
\providecommand{\href}[2]{#2}

\appendix
\section{Generalised toric codes over $\F_q$, $q\leq 7$}\label{app:toric_table}
The largest minimum distance $d_t$ achieved by a generalised toric code with given block length $n$ and dimension $k$ over the field $\F_q$ are documented in the following table. For comparison, we also give the current best minimum distance $d_g$ amongst all linear codes as given by Grassl~\cite{grassl}. The penultimate column records whether $d_t\ge d_g$, and the final column gives an example collection of points realising these invariants. The only example for which a \emph{generalised} toric code is required is $[36,7,23]$ over $\F_7$.

\begin{center}
\begin{longtable}{ccccccp{4.25in}}
\toprule
$q$&$n$&$k$&$d_t$&$d_g$&Best&\multicolumn{1}{c}{Example points}\\
\cmidrule(lr){1-5}\cmidrule(lr){6-6}\cmidrule(lr){7-7}
\endfirsthead
\multicolumn{7}{l}{\small\emph{continued from previous page}}\\
\toprule
$q$&$n$&$k$&$d_t$&$d_g$&Best&\multicolumn{1}{c}{Example points}\\
\cmidrule(lr){1-5}\cmidrule(lr){6-6}\cmidrule(lr){7-7}
\endhead
\bottomrule
\multicolumn{7}{r}{\small\emph{continued on next page}}\\
\endfoot
\endlastfoot
$3$&$4$&$3$&$2$&$2$&y&\small$(0,1),$ $(1,0),$ $(1,1)$\\
$3$&$4$&$4$&$1$&$1$&y&\small$(0,0),$ $(0,1),$ $(1,0),$ $(1,1)$\\
\cmidrule(lr){1-5}\cmidrule(lr){6-6}\cmidrule(lr){7-7}
$4$&$9$&$3$&$6$&$6$&y&\small$(0,1),$ $(1,0),$ $(1,1)$\\
$4$&$9$&$4$&$4$&$5$&n&\small$(0,0),$ $(0,1),$ $(1,0),$ $(1,1)$\\
$4$&$9$&$5$&$3$&$4$&n&\small$(0,1),$ $(0,2),$ $(1,1),$ $(1,2),$ $(2,0)$\\
$4$&$9$&$6$&$3$&$3$&y&\small$(0,1),$ $(0,2),$ $(1,1),$ $(1,2),$ $(2,0),$ $(2,1)$\\
$4$&$9$&$7$&$2$&$2$&y&\small$(0,1),$ $(0,2),$ $(1,1),$ $(1,2),$ $(2,0),$ $(2,1),$ $(2,2)$\\
$4$&$9$&$8$&$2$&$2$&y&\small$(0,1),$ $(0,2),$ $(1,0),$ $(1,1),$ $(1,2),$ $(2,0),$ $(2,1),$ $(2,2)$\\
$4$&$9$&$9$&$1$&$1$&y&\small$(0,0),$ $(0,1),$ $(0,2),$ $(1,0),$ $(1,1),$ $(1,2),$ $(2,0),$ $(2,1),$ $(2,2)$\\
\cmidrule(lr){1-5}\cmidrule(lr){6-6}\cmidrule(lr){7-7}
$5$&$16$&$3$&$12$&$12$&y&\small$(0,1),$ $(1,0),$ $(1,1)$\\
$5$&$16$&$4$&$10$&$11$&n&\small$(0,1),$ $(1,1),$ $(1,2),$ $(2,0)$\\
$5$&$16$&$5$&$8$&$9$&n&\small$(0,1),$ $(0,2),$ $(1,1),$ $(1,2),$ $(2,0)$\\
$5$&$16$&$6$&$8$&$8$&y&\small$(0,1),$ $(0,2),$ $(1,1),$ $(1,2),$ $(2,0),$ $(2,1)$\\
$5$&$16$&$7$&$7$&$7$&y&\small$(0,3),$ $(1,2),$ $(1,3),$ $(2,0),$ $(2,1),$ $(2,2),$ $(3,2)$\\
$5$&$16$&$8$&$6$&$7$&n&\small$(0,1),$ $(0,2),$ $(1,0),$ $(1,1),$ $(1,2),$ $(2,0),$ $(2,1),$ $(2,2)$\\
$5$&$16$&$9$&$6$&$6$&y&\small$(0,3),$ $(1,1),$ $(1,2),$ $(1,3),$ $(2,0),$ $(2,1),$ $(2,2),$ $(3,1),$ $(3,2)$\\
$5$&$16$&$10$&$4$&$5$&n&\small$(0,3),$ $(1,0),$ $(1,1),$ $(1,2),$ $(1,3),$ $(2,1),$ $(2,2),$ $(2,3),$ $(3,1),$ $(3,2)$\\
$5$&$16$&$11$&$4$&$4$&y&\small$(0,3),$ $(1,0),$ $(1,1),$ $(1,2),$ $(1,3),$ $(2,1),$ $(2,2),$ $(2,3),$ $(3,1),$ $(3,2),$ $(3,3)$\\
$5$&$16$&$12$&$4$&$4$&y&\small$(0,2),$ $(0,3),$ $(1,1),$ $(1,2),$ $(1,3),$ $(2,0),$ $(2,1),$ $(2,2),$ $(2,3),$ $(3,0),$ $(3,1),$ $(3,2)$\\
$5$&$16$&$13$&$3$&$3$&y&\small$(0,1),$ $(0,2),$ $(1,0),$ $(1,1),$ $(1,2),$ $(1,3),$ $(2,0),$ $(2,1),$ $(2,2),$ $(2,3),$ $(3,0),$ $(3,1),$ $(3,2)$\\
$5$&$16$&$14$&$2$&$2$&y&\small$(0,1),$ $(0,2),$ $(1,0),$ $(1,1),$ $(1,2),$ $(1,3),$ $(2,0),$ $(2,1),$ $(2,2),$ $(2,3),$ $(3,0),$ $(3,1),$ $(3,2),$ $(3,3)$\\
$5$&$16$&$15$&$2$&$2$&y&\small$(0,1),$ $(0,2),$ $(0,3),$ $(1,0),$ $(1,1),$ $(1,2),$ $(1,3),$ $(2,0),$ $(2,1),$ $(2,2),$ $(2,3),$ $(3,0),$ $(3,1),$ $(3,2),$ $(3,3)$\\
$5$&$16$&$16$&$1$&$1$&y&\small$(0,0),$ $(0,1),$ $(0,2),$ $(0,3),$ $(1,0),$ $(1,1),$ $(1,2),$ $(1,3),$ $(2,0),$ $(2,1),$ $(2,2),$ $(2,3),$ $(3,0),$ $(3,1),$ $(3,2),$ $(3,3)$\\
\cmidrule(lr){1-5}\cmidrule(lr){6-6}\cmidrule(lr){7-7}
$7$&$36$&$3$&$30$&$30$&y&\small$(0,1),$ $(1,0),$ $(1,1)$\\
$7$&$36$&$4$&$27$&$28$&n&\small$(0,1),$ $(1,1),$ $(1,2),$ $(2,0)$\\
$7$&$36$&$5$&$24$&$27$&n&\small$(0,1),$ $(0,2),$ $(1,1),$ $(1,2),$ $(2,0)$\\
$7$&$36$&$6$&$24$&$25$&n&\small$(0,1),$ $(0,2),$ $(1,1),$ $(1,2),$ $(2,0),$ $(2,1)$\\
$7$&$36$&$7$&$23$&$24$&n&\small$(0,3),$ $(1,2),$ $(1,3),$ $(2,0),$ $(2,2),$ $(3,1),$ $(3,2)$\\
$7$&$36$&$8$&$20$&$22$&n&\small$(0,1),$ $(0,2),$ $(1,0),$ $(1,1),$ $(1,2),$ $(2,0),$ $(2,1),$ $(2,2)$\\
$7$&$36$&$9$&$20$&$21$&n&\small$(0,3),$ $(1,1),$ $(1,2),$ $(1,3),$ $(2,0),$ $(2,1),$ $(2,2),$ $(3,1),$ $(3,2)$\\
$7$&$36$&$10$&$18$&$20$&n&\small$(0,3),$ $(1,0),$ $(1,1),$ $(1,2),$ $(1,3),$ $(2,1),$ $(2,2),$ $(2,3),$ $(3,1),$ $(3,2)$\\
$7$&$36$&$11$&$18$&$19$&n&\small$(0,4),$ $(1,3),$ $(1,4),$ $(2,1),$ $(2,2),$ $(2,3),$ $(3,0),$ $(3,1),$ $(3,2),$ $(3,3),$ $(4,3)$\\
$7$&$36$&$12$&$17$&$18$&n&\small$(0,4),$ $(1,3),$ $(1,4),$ $(2,1),$ $(2,2),$ $(2,3),$ $(3,0),$ $(3,1),$ $(3,2),$ $(3,3),$ $(4,2),$ $(4,3)$\\
$7$&$36$&$13$&$15$&$17$&n&\small$(0,1),$ $(0,2),$ $(1,0),$ $(1,1),$ $(1,2),$ $(1,3),$ $(2,0),$ $(2,1),$ $(2,2),$ $(2,3),$ $(3,0),$ $(3,1),$ $(3,2)$\\
$7$&$36$&$14$&$15$&$16$&n&\small$(1,1),$ $(1,2),$ $(1,3),$ $(1,4),$ $(2,1),$ $(2,2),$ $(2,3),$ $(2,4),$ $(3,1),$ $(3,2),$ $(3,3),$ $(4,0),$ $(4,1),$ $(4,2)$\\
$7$&$36$&$15$&$14$&$15$&n&\small$(0,4),$ $(1,2),$ $(1,3),$ $(1,4),$ $(2,1),$ $(2,2),$ $(2,3),$ $(2,4),$ $(3,0),$ $(3,1),$ $(3,2),$ $(3,3),$ $(4,1),$ $(4,2),$ $(4,3)$\\
$7$&$36$&$16$&$12$&$14$&n&\small$(0,3),$ $(1,0),$ $(1,1),$ $(1,2),$ $(1,3),$ $(2,1),$ $(2,2),$ $(2,3),$ $(3,1),$ $(3,2),$ $(3,3),$ $(3,4),$ $(4,1),$ $(4,2),$ $(4,3),$ $(4,4)$\\
$7$&$36$&$17$&$12$&$13$&n&\small$(0,3),$ $(1,0),$ $(1,1),$ $(1,2),$ $(1,3),$ $(1,4),$ $(2,1),$ $(2,2),$ $(2,3),$ $(2,4),$ $(3,1),$ $(3,2),$ $(3,3),$ $(3,4),$ $(4,2),$ $(4,3),$ $(4,4)$\\
$7$&$36$&$18$&$12$&$13$&n&\small$(0,3),$ $(0,4),$ $(1,2),$ $(1,3),$ $(1,4),$ $(2,0),$ $(2,1),$ $(2,2),$ $(2,3),$ $(2,4),$ $(3,0),$ $(3,1),$ $(3,2),$ $(3,3),$ $(4,0),$ $(4,1),$ $(4,2),$ $(4,3)$\\
$7$&$36$&$19$&$12$&$11$&y&\small$(0,4),$ $(1,3),$ $(1,4),$ $(2,2),$ $(2,3),$ $(2,4),$ $(2,5),$ $(3,2),$ $(3,3),$ $(3,4),$ $(3,5),$ $(4,1),$ $(4,2),$ $(4,3),$ $(4,4),$ $(4,5),$ $(5,0),$ $(5,1),$ $(5,2)$\\
$7$&$36$&$20$&$10$&$11$&n&\small$(0,3),$ $(1,2),$ $(1,3),$ $(2,1),$ $(2,2),$ $(2,3),$ $(2,4),$ $(3,1),$ $(3,2),$ $(3,3),$ $(3,4),$ $(3,5),$ $(4,1),$ $(4,2),$ $(4,3),$ $(4,4),$ $(4,5),$ $(5,0),$ $(5,1),$ $(5,2)$\\
$7$&$36$&$21$&$10$&$10$&y&\small$(0,1),$ $(0,2),$ $(0,3),$ $(0,4),$ $(1,1),$ $(1,2),$ $(1,3),$ $(1,4),$ $(1,5),$ $(2,1),$ $(2,2),$ $(2,3),$ $(2,4),$ $(2,5),$ $(3,1),$ $(3,2),$ $(3,3),$ $(4,0),$ $(4,1),$ $(4,2),$ $(5,0)$\\
$7$&$36$&$22$&$9$&$10$&n&\small$(0,3),$ $(0,4),$ $(1,2),$ $(1,3),$ $(1,4),$ $(2,1),$ $(2,2),$ $(2,3),$ $(2,4),$ $(2,5),$ $(3,1),$ $(3,2),$ $(3,3),$ $(3,4),$ $(3,5),$ $(4,0),$ $(4,1),$ $(4,2),$ $(4,3),$ $(5,0),$ $(5,1),$ $(5,2)$\\
$7$&$36$&$23$&$8$&$9$&n&\small$(1,2),$ $(1,3),$ $(1,4),$ $(1,5),$ $(2,1),$ $(2,2),$ $(2,3),$ $(2,4),$ $(2,5),$ $(3,1),$ $(3,2),$ $(3,3),$ $(3,4),$ $(3,5),$ $(4,0),$ $(4,1),$ $(4,2),$ $(4,3),$ $(4,4),$ $(5,1),$ $(5,2),$ $(5,3),$ $(5,4)$\\
$7$&$36$&$24$&$8$&$8$&y&\small$(0,3),$ $(0,4),$ $(1,2),$ $(1,3),$ $(1,4),$ $(1,5),$ $(2,1),$ $(2,2),$ $(2,3),$ $(2,4),$ $(2,5),$ $(3,1),$ $(3,2),$ $(3,3),$ $(3,4),$ $(3,5),$ $(4,0),$ $(4,1),$ $(4,2),$ $(4,3),$ $(4,4),$ $(5,0),$ $(5,1),$ $(5,2)$\\
$7$&$36$&$25$&$6$&$7$&n&\small$(0,1),$ $(0,2),$ $(0,3),$ $(1,1),$ $(1,2),$ $(1,3),$ $(1,4),$ $(1,5),$ $(2,1),$ $(2,2),$ $(2,3),$ $(2,4),$ $(2,5),$ $(3,1),$ $(3,2),$ $(3,3),$ $(3,4),$ $(3,5),$ $(4,0),$ $(4,1),$ $(4,2),$ $(4,3),$ $(4,4),$ $(4,5),$ $(5,4)$\\
$7$&$36$&$26$&$6$&$7$&n&\small$(0,1),$ $(0,2),$ $(0,3),$ $(0,4),$ $(1,1),$ $(1,2),$ $(1,3),$ $(1,4),$ $(1,5),$ $(2,1),$ $(2,2),$ $(2,3),$ $(2,4),$ $(2,5),$ $(3,1),$ $(3,2),$ $(3,3),$ $(3,4),$ $(3,5),$ $(4,0),$ $(4,1),$ $(4,2),$ $(4,3),$ $(4,4),$ $(4,5),$ $(5,4)$\\
$7$&$36$&$27$&$6$&$6$&y&\small$(0,3),$ $(0,4),$ $(0,5),$ $(1,2),$ $(1,3),$ $(1,4),$ $(1,5),$ $(2,1),$ $(2,2),$ $(2,3),$ $(2,4),$ $(2,5),$ $(3,0),$ $(3,1),$ $(3,2),$ $(3,3),$ $(3,4),$ $(3,5),$ $(4,0),$ $(4,1),$ $(4,2),$ $(4,3),$ $(4,4),$ $(5,0),$ $(5,1),$ $(5,2),$ $(5,3)$\\
$7$&$36$&$28$&$6$&$6$&y&\small$(0,3),$ $(0,4),$ $(0,5),$ $(1,1),$ $(1,2),$ $(1,3),$ $(1,4),$ $(1,5),$ $(2,1),$ $(2,2),$ $(2,3),$ $(2,4),$ $(2,5),$ $(3,0),$ $(3,1),$ $(3,2),$ $(3,3),$ $(3,4),$ $(3,5),$ $(4,0),$ $(4,1),$ $(4,2),$ $(4,3),$ $(4,4),$ $(5,0),$ $(5,1),$ $(5,2),$ $(5,3)$\\
$7$&$36$&$29$&$5$&$6$&n&\small$(0,2),$ $(0,3),$ $(0,4),$ $(0,5),$ $(1,2),$ $(1,3),$ $(1,4),$ $(1,5),$ $(2,1),$ $(2,2),$ $(2,3),$ $(2,4),$ $(2,5),$ $(3,0),$ $(3,1),$ $(3,2),$ $(3,3),$ $(3,4),$ $(3,5),$ $(4,0),$ $(4,1),$ $(4,2),$ $(4,3),$ $(4,4),$ $(4,5),$ $(5,0),$ $(5,1),$ $(5,2),$ $(5,3)$\\
$7$&$36$&$30$&$4$&$5$&n&\small$(0,1),$ $(0,2),$ $(0,3),$ $(0,4),$ $(0,5),$ $(1,1),$ $(1,2),$ $(1,3),$ $(1,4),$ $(1,5),$ $(2,0),$ $(2,1),$ $(2,2),$ $(2,3),$ $(2,4),$ $(2,5),$ $(3,0),$ $(3,1),$ $(3,2),$ $(3,3),$ $(3,4),$ $(3,5),$ $(4,1),$ $(4,2),$ $(4,3),$ $(4,4),$ $(4,5),$ $(5,2),$ $(5,3),$ $(5,4)$\\
$7$&$36$&$31$&$4$&$4$&y&\small$(0,1),$ $(0,2),$ $(0,3),$ $(0,4),$ $(0,5),$ $(1,1),$ $(1,2),$ $(1,3),$ $(1,4),$ $(1,5),$ $(2,1),$ $(2,2),$ $(2,3),$ $(2,4),$ $(2,5),$ $(3,0),$ $(3,1),$ $(3,2),$ $(3,3),$ $(3,4),$ $(3,5),$ $(4,0),$ $(4,1),$ $(4,2),$ $(4,3),$ $(4,4),$ $(4,5),$ $(5,0),$ $(5,1),$ $(5,2),$ $(5,3)$\\
$7$&$36$&$32$&$3$&$4$&n&\small$(0,1),$ $(0,2),$ $(0,3),$ $(0,4),$ $(0,5),$ $(1,1),$ $(1,2),$ $(1,3),$ $(1,4),$ $(1,5),$ $(2,0),$ $(2,1),$ $(2,2),$ $(2,3),$ $(2,4),$ $(2,5),$ $(3,0),$ $(3,1),$ $(3,2),$ $(3,3),$ $(3,4),$ $(3,5),$ $(4,0),$ $(4,1),$ $(4,2),$ $(4,3),$ $(4,4),$ $(4,5),$ $(5,1),$ $(5,2),$ $(5,3),$ $(5,4)$\\
$7$&$36$&$33$&$3$&$3$&y&\small$(0,1),$ $(0,2),$ $(0,3),$ $(0,4),$ $(0,5),$ $(1,1),$ $(1,2),$ $(1,3),$ $(1,4),$ $(1,5),$ $(2,0),$ $(2,1),$ $(2,2),$ $(2,3),$ $(2,4),$ $(2,5),$ $(3,0),$ $(3,1),$ $(3,2),$ $(3,3),$ $(3,4),$ $(3,5),$ $(4,0),$ $(4,1),$ $(4,2),$ $(4,3),$ $(4,4),$ $(4,5),$ $(5,0),$ $(5,1),$ $(5,2),$ $(5,3),$ $(5,4)$\\
$7$&$36$&$34$&$2$&$2$&y&\small$(0,1),$ $(0,2),$ $(0,3),$ $(0,4),$ $(0,5),$ $(1,0),$ $(1,1),$ $(1,2),$ $(1,3),$ $(1,4),$ $(1,5),$ $(2,0),$ $(2,1),$ $(2,2),$ $(2,3),$ $(2,4),$ $(2,5),$ $(3,0),$ $(3,1),$ $(3,2),$ $(3,3),$ $(3,4),$ $(3,5),$ $(4,0),$ $(4,1),$ $(4,2),$ $(4,3),$ $(4,4),$ $(4,5),$ $(5,1),$ $(5,2),$ $(5,3),$ $(5,4),$ $(5,5)$\\
$7$&$36$&$35$&$2$&$2$&y&\small$(0,1),$ $(0,2),$ $(0,3),$ $(0,4),$ $(0,5),$ $(1,0),$ $(1,1),$ $(1,2),$ $(1,3),$ $(1,4),$ $(1,5),$ $(2,0),$ $(2,1),$ $(2,2),$ $(2,3),$ $(2,4),$ $(2,5),$ $(3,0),$ $(3,1),$ $(3,2),$ $(3,3),$ $(3,4),$ $(3,5),$ $(4,0),$ $(4,1),$ $(4,2),$ $(4,3),$ $(4,4),$ $(4,5),$ $(5,0),$ $(5,1),$ $(5,2),$ $(5,3),$ $(5,4),$ $(5,5)$\\
$7$&$36$&$36$&$1$&$1$&y&\small$(0,0),$ $(0,1),$ $(0,2),$ $(0,3),$ $(0,4),$ $(0,5),$ $(1,0),$ $(1,1),$ $(1,2),$ $(1,3),$ $(1,4),$ $(1,5),$ $(2,0),$ $(2,1),$ $(2,2),$ $(2,3),$ $(2,4),$ $(2,5),$ $(3,0),$ $(3,1),$ $(3,2),$ $(3,3),$ $(3,4),$ $(3,5),$ $(4,0),$ $(4,1),$ $(4,2),$ $(4,3),$ $(4,4),$ $(4,5),$ $(5,0),$ $(5,1),$ $(5,2),$ $(5,3),$ $(5,4),$ $(5,5)$\\
\bottomrule
\end{longtable}
\end{center}
\end{document}